\documentclass{amsart}

\usepackage{amsthm}
\usepackage{amsmath}
\usepackage{amsfonts}
\usepackage{graphicx}

\newcommand{\la}{\lambda}
\newcommand{\N}{\mathbb{N}}
\providecommand{\abs}[1]{\lvert#1\rvert}%

\newtheorem{theorem}{Theorem}
\newtheorem{conjecture}[theorem]{Conjecture}
\theoremstyle{definition}
\newtheorem{definition}[theorem]{Definition}
\newtheorem{observation}[theorem]{Observation}

\usepackage[pdftitle={Strict Partitions of Maximal Projective Degree},pdfauthor={Dan Bernstein},bookmarks=true]{hyperref}

\begin{document}

\title{Strict Partitions of Maximal Projective Degree}
\date{May 28, 2007}
\author{Dan Bernstein}
\address{
    Department of Mathematics\\
    The Weizmann Institute of Science\\
    Rehovot 76100, Israel
}
\email{dan.bernstein@weizmann.ac.il}
\thanks{Partially supported by the Israel Science Foundation (grant No.~947/04).}

\begin{abstract}
The projective degrees of strict partitions of $n$ were computed for all $n\le 100$ and
the partitions with maximal projective degree were found for each $n$. It was observed that maximizing
partitions for successive values of $n$ ``lie close to each other'' in a certain sense.
Conjecturing that this holds for larger values of $n$, the partitions of maximal degree
were computed for all $n \le 220$. The results are consistent with a recent conjecture on
the limiting shape of the strict partition of maximal projective degree.
\end{abstract}

\maketitle

\section{Introduction}

Let $\la=(\la_1,\,\la_2,\,\dots)$ be a partition of $n$, denoted as usual by $\la \vdash n$. Let $f^\la$ denote the number of standard tableaux of shape $\la$. $f^\la$ is also the number of paths in the Young graph $Y$ from the root $(1)$ to $\la$, and it is also the degree of the irreducible character $\chi^\la$ of the symmetric group $S_n$.

A partition $\la=(\la_1,\,\dots,\,\la_r)\vdash n$ is said to be \emph{strict} if
$\la_1 > \la_2 > \cdots > \la_r > 0$ for some $r$. In that case we write $\la \models n$.
Let $SY$ be the subgraph of the Young graph $Y$ formed by the strict partitions. If $\la\models n$, let $g^\la$ denote the number of paths in $SY$ from the root $(1)$ to $\la$. According to a theorem by Schur, $g^\la$ is the degree of the projective representation of $S_n$ corresponding to $\la$. It is also the number of standard young tableaux of \emph{shifted shape} $\la$.

Vershik and Kerov~\cite{vershik:77, vershik:85} have determined the asymptotic shape of the partition $\la$ that maximizes $f^\la$ as $\abs{\la} = n\to\infty$. They have also shown  that the same shape is also the asymptotic expected shape of a random partition with respect to the Plancherel measure. The latter result was reached independently by Logan and Shepp~\cite{logan:variational} as well. Through the Robinson-Schensted algorithm, the expected shape relates to the expected length of the longest increasing subsequence in a random permutation. For some recent developments related to this problem and the probability distributions involved, see~\cite{baik:longest, baik:second, stanley:recent}

More precisely, given the Young diagram of a partition $\la \vdash n$ where each box is $1\times 1$, shrink it along both axes by a factor of $\sqrt{n}$ to obtain the re-scaled diagram $\bar\la$ of total area $1$. For each $n$, let $\la^{(n)}_{f\max}$ be a partition $\la\vdash n$ with maximal $f^\la$, that is, $f^{\la^{(n)}_{f\max}} = \max\{\;f^\nu \mid \nu \vdash n\;\}$. Through slight abuse of notation, where the maximizing partition is not unique for a given $n$, we shall take $\la^{(n)}_{f\max}$ to read ``any $\la\vdash n$ of maximal degree''.

\begin{theorem}[\cite{vershik:77, vershik:85}]\label{TH:VK}
The limit shape as $n\to\infty$ of the re-scaled diagrams $\bar\la^{(n)}_{f\max}$ exists, and is given by the two axes and by the parametric curve
\[
\begin{cases}
	x=\left(\frac{2}{\pi}\right)\left( \sin\Theta - \Theta \cos \Theta \right) + 2 \cos \Theta,\\
	y = -\left( \frac{2}{\pi} \right) \left( \sin\Theta - \Theta \cos\Theta \right),
\end{cases}
\quad 0 \le \Theta \le \pi	.
\]
\end{theorem}

\begin{figure}[htb]
\begin{center}
\includegraphics[width=1.25in, viewport=72 0 216 288]{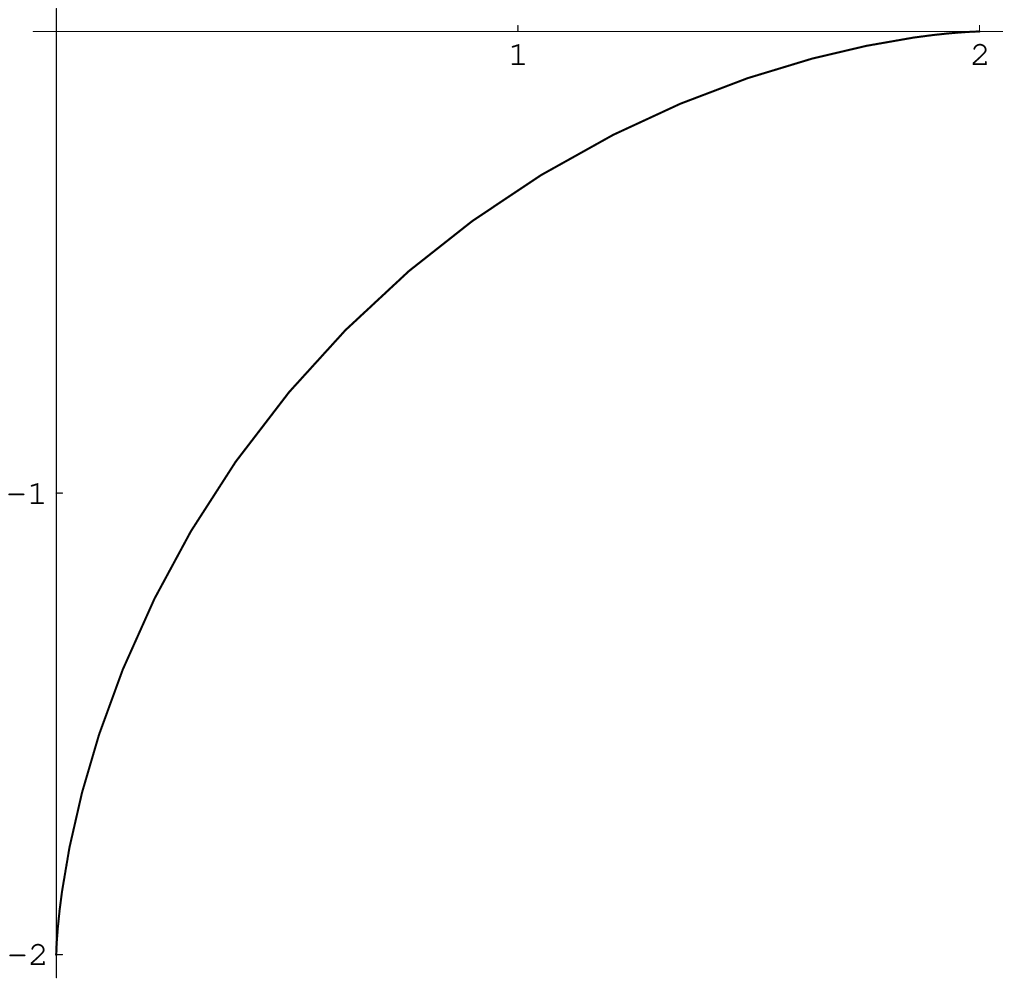}
\end{center}
\refstepcounter{figure}
\centerline{Figure~\arabic{figure}}
\label{FIG:VK}
\end{figure}

The Vershik-Kerov limit shape of Theorem~\ref{TH:VK} is shown in Figure~\ref{FIG:VK}.

The problem of determining the asymptotic shape of the partition $\la$ which maximizes $g^\la$ remains unsolved, and we are unaware of even partial characterizations of the shape. However, recently the following was conjectured.

\begin{conjecture}[{\cite[Conjecture~8.2]{henke:maximal}}]\label{CONJ:BHR}
The limit shape $\la^*$ of the $\la \models n$ maximizing $2^{n-\ell(\la)}\left(g^\la\right)^2$ --- and possibly maximizing $g^\la$ --- is given by
the two axes and by the parametric curve
\[
\begin{cases}
	x = 2\sqrt{2}\cos\Theta,\\
	y=\left(\frac{2\sqrt{2}}{\pi}\right)\left(\Theta \cos\Theta - \sin\Theta\right),
\end{cases}
0 \le \Theta \le \frac{\pi}{2}	.
\]
\end{conjecture}

\begin{figure}[htb]
\begin{center}
\includegraphics[width=1.75in, viewport=72 0 216 108]{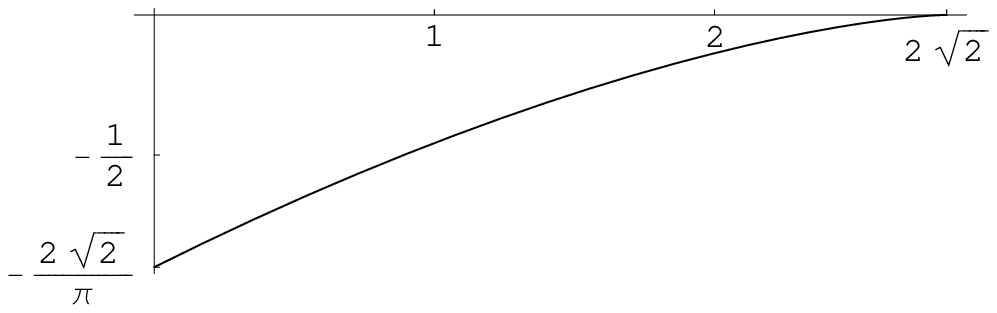}
\end{center}
\refstepcounter{figure}
\centerline{Figure~\arabic{figure}}
\label{FIG:BHR}
\end{figure}

The conjectured limit shape of Conjecture~\ref{CONJ:BHR} is shown in Figure~\ref{FIG:BHR}. It was
obtained from the Vershik-Kerov shape by bisecting it along the line $y=-x$, taking the upper half, dilating it (so that its area become $1$) and applying the shearing transformation $(x,y)\mapsto(x+y,y)$ (to bring the line $y=-x$ to the $y$ axis).

In the next section, we give the results of computing the partition maximizing $g^\la$ over all $\la\models n$ for $1 \le n \le 100$. We observe a property of successive maximizing partitions in the range $1 \le n \le 100$, and conjecture that it holds for all $n$. Assuming the conjecture, we compute the maximizing partitions for $100 < n \le 220$. Our results are consistent with Conjecture~\ref{CONJ:BHR} for both $g^\la$ and $2^{-\ell(\la)}\left(g^\la\right)^2$.

\section{Results for $n \le 100$}

To compute $g^\la$, we used the following formula, due to Schur.
\begin{theorem}
Let $\la=(\la_1,\,\dots,\,\la_r)\models n$. Then
\[
	g^\la = \frac{n!}{\la_1!\cdots\la_r!}\prod_{1 \le i < j \le r}\frac{\la_i-\la_j}{\la_i+\la_j}
\]
\end{theorem}

All computations were done with {\it Mathematica}.

For $1 \le n \le 100$, $g^\la$ was computed for every $\la \models n$ and the partitions attaining the maximal value for each $n$ were identified.

\begin{observation}[Uniqueness of the maximum]
For every $n\in [100]\setminus\{3,11\}$, there exists a partition $\la=\la^{(n)}_{g\max}\models n$
such that $g^\la > g^\mu$ for all $\la \neq \mu \models n$.
\end{observation}

\begin{observation}
\[
\la^{(100)}_{g\max} = (24,\, 20,\, 16,\, 13,\, 10,\, 8,\, 5,\, 3,\, 1)	.
\]
\end{observation}

\begin{figure}[htb]
\begin{center}
\includegraphics[width=1.75in, viewport=72 0 216 108]{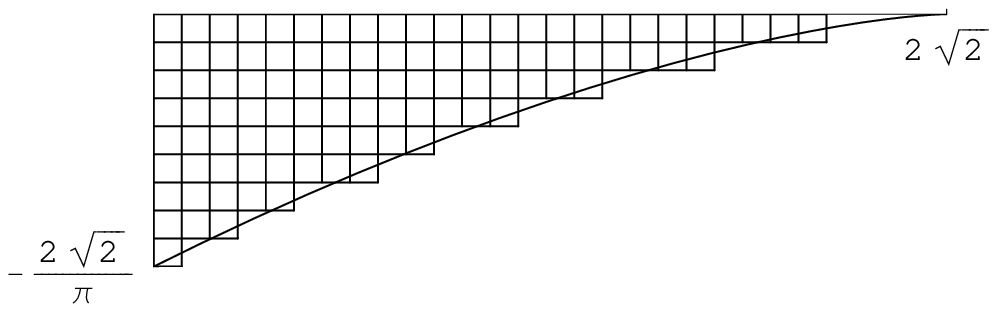}
\end{center}
\refstepcounter{figure}
\centerline{Figure~\arabic{figure}}
\label{FIG:100}
\end{figure}

Figure~\ref{FIG:100} shows the normalized diagram $\bar\la^{(100)}_{g\max}$ overlaid with the conjectured limit shape.

Recall that by Conjecture~\ref{CONJ:BHR}, the partitions maximizing $g^\la$ are asymptotically equal to the partitions maximizing $2^{-\ell(\la)}\left(g^\la\right)^2$. For all strict partitions $\la$ of $1 \le n \le 100$, $2^{-\ell(\la)}\left(g^\la\right)^2$ was computed and the partition maximizing it was denoted $\la^{(n)}_{2g\max}$. Comparing the results with the results for $g^\la$, the following was observed. 
\begin{observation}
For all
\[
n \in \{1,\,2,\,\dots,\,100\} \setminus \{3,\, 8,\, 16,\, 25,\, 26,\, 38,\, 51,\, 52,\, 53,\, 54,\, 69,\, 70,\, 88,\, 89,\, 90,\, 91\}	,
\]
 $\la^{(n)}_{g\max} = \la^{(n)}_{2g\max}$.
\end{observation}

\begin{definition}[$d$-successor]
Let $\mu  = (\mu_1,\,\dots,\,\mu_r) \models n$ and let $d\in\N$. The $d$-successors of $\mu$
are the elements of the set
\[
	N(\mu,d) := \{\; \la=(\la_1,\,\dots,\,\la_s)\models n+1 \mid \abs{\la_i-\mu_i}\le d\quad 1\le i \le s \;\}
\]
\end{definition}

We have observed the following.
\begin{observation}
For $1 \le n < 100$, $\la^{(n+1)}_{g\max}\in N(\la^{n}_{g\max},1)$.
\end{observation}
When $\la^{(n)}_{g\max}$ is not unique, read the above to mean ``every $\la\models n+1$ of maximal projective degree is a $1$-successor of every $\la\models n$ of maximal degree''.

\section{A Conjecture and Results for $n \le 220$}

Based on the above observation, we conjecture the following.
\begin{conjecture}[Maximizers are successors to a maximizers]\label{CONJ:succ}
For all $n$, if $\la^{(n+1)}_{g\max}\in N(\la^{n}_{g\max},1)$. 
\end{conjecture}

Assuming that the conjecture holds, $\la^{(n+1)}_{g\max}$ was computed for $100 \le n < 220$ as follows: statring with $\la=\la^{(n)}_{g\max}$, for every $\mu \in N(\la,1)$, the ratio
$\frac{g^\mu}{g^\la}$ was computed and the $\mu$ maximizing the ratio was selected.

\begin{observation}
If Conjecture~\ref{CONJ:succ} holds for all $n < 220$, then
\[
\la^{(220)}_{g\max} = (37,\, 32,\, 28,\, 24,\, 21,\, 18,\, 16,\, 13,\, 11,\, 8,\, 6,\, 4,\, 2)	.
\]
\end{observation}

\begin{figure}[h]
\begin{center}
\includegraphics[width=1.75in, viewport=72 0 216 108]{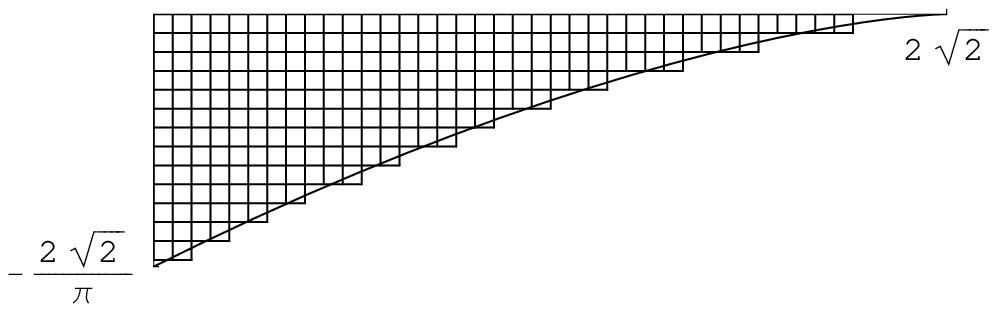}
\end{center}
\refstepcounter{figure}
\centerline{Figure~\arabic{figure}}
\label{FIG:220}
\end{figure}

Figure~\ref{FIG:220} shows the normalized diagram $\bar\la^{(220)}_{g\max}$ overlaid with the conjectured limit shape.

\begin{figure}[h]
\begin{center}
\includegraphics[width=2in, viewport=72 0 216 216]{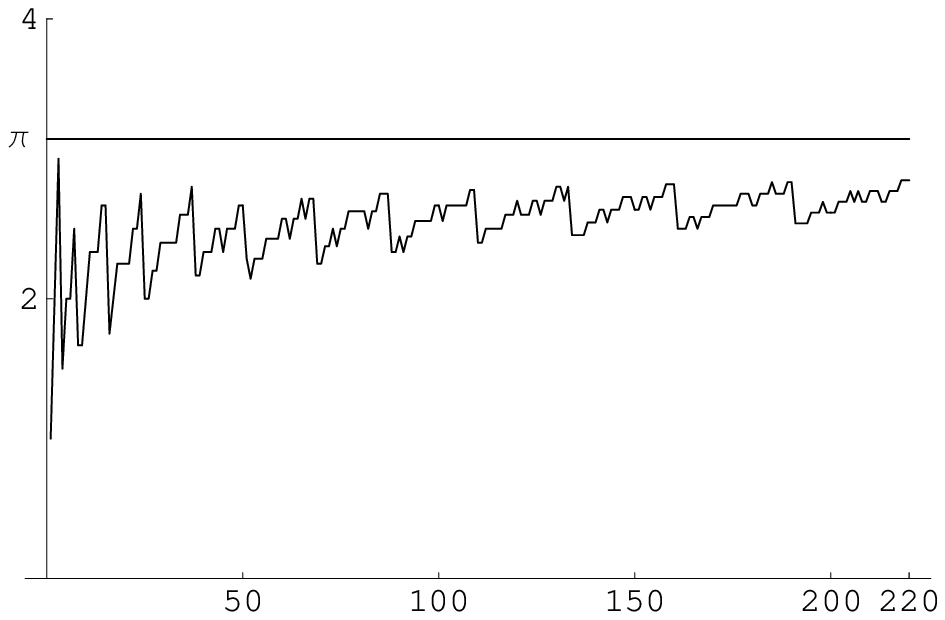}
\end{center}
\refstepcounter{figure}
\centerline{Figure~\arabic{figure}}
\label{FIG:ratio}
\end{figure}

Figure~\ref{FIG:ratio} shows the ratio $\frac{(\la^{(n)}_{g\max})_1}{\ell(\la^{(n)}_{g\max})}$ for $1 \le n \le 220$, where the values above $n=100$ are based on Conjecture~\ref{CONJ:succ}. According to Conjecture~\ref{CONJ:BHR}, the ratio at the limit is $\pi$.



\newpage
\begin{thebibliography}{1}

\bibitem{baik:longest}
J.~Baik, P.~Deift, and K.~Johansson.
\newblock On the distribution of the length of the longest increasing
  subsequence of random permutations.
\newblock {\em J. Amer. Math. Soc.}, 12(4):1119--1178, 1999.

\bibitem{baik:second}
J.~Baik, P.~Deift, and K.~Johansson.
\newblock On the distribution of the length of the second row of a {Y}oung
  diagram under {P}lancherel measure.
\newblock {\em Geom. Funct. Anal.}, 10(4):702--731, 2000.

\bibitem{henke:maximal}
D.~Bernstein, A.~Henke, and A.~Regev.
\newblock Maximal projective degrees for strict partitions.
\newblock Preprint, 2006.

\bibitem{logan:variational}
B.~F. Logan and L.~A. Shepp.
\newblock A variational problem for random {Y}oung tableaux.
\newblock {\em Advances in Math.}, 26(2):206--222, 1977.

\bibitem{stanley:recent}
R.~P. Stanley.
\newblock Recent progress in algebraic combinatorics.
\newblock {\em Bull. Amer. Math. Soc. (N.S.)}, 40(1):55--68 (electronic), 2003.
\newblock Mathematical challenges of the 21st century (Los Angeles, CA, 2000).

\bibitem{vershik:77}
A.~M. Vershik and S.~V. Kerov.
\newblock Asymptotics of the {P}lancherel measure of the symmetric group and
  the limit form of young tableaux.
\newblock {\em Dokl. Akad. Nauk SSSR}, 233:1024--1027, 1977.

\bibitem{vershik:85}
A.~M. Vershik and S.~V. Kerov.
\newblock Asymptotic behavior of the maximum and generic dimensions of
  irreducible representations of the symmetric group.
\newblock {\em Funktsional. Anal. i Prilozhen.}, 19(1):25--36, 96, 1985.

\end{thebibliography}
\end{document}